%% file: article_mirekfrost_martinkruzik_janvaldman.tex
\documentclass[a4paper,11pt]{article}
\pdfoutput=1

\usepackage{color}   
\usepackage{amsmath,amsfonts,amssymb,amsthm,authblk}
\usepackage{graphicx}
\usepackage{epsfig}
\usepackage{multicol}

\usepackage{hyperref}
\hypersetup{colorlinks=true, breaklinks=true, linkcolor=black, menucolor=black, urlcolor=black, citecolor=black}

\usepackage{soul}
\newtheorem{theorem}{Theorem}[section]
\newtheorem{definition}[theorem]{Definition}

\newtheorem{remark}[theorem]{Remark}

\usepackage{todonotes}
\usepackage{comment}

\definecolor{darkgreen}{rgb}{0,0.5,0}

\input{all-def}

\newcommand{\quadlown}{\qquad \qquad}

 \newcommand{\md}{{\rm \mbox{d}}}
 \renewcommand{\O}{\Omega}
 \newcommand{\be}{\begin{eqnarray}}
 \newcommand{\ee}{\end{eqnarray}}

\DeclareMathOperator{\tr}{tr}

\title{Interfacial polyconvex energy-enhanced evolutionary model for shape memory alloys}
\author[1]{Miroslav Frost}
\author[2,3]{Martin Kru\v{z}\'{\i}k}
\author[2,4]{Jan Valdman}
\affil[1]{Czech Academy of Sciences, Institute of Thermomechanics, Dolej\v{s}kova 5, CZ-182 00 Prague 8, Czechia,}
\affil[2]{Czech academy of Sciences, Institute of Information Theory and Automation,  Pod vod\'{a}renskou v\v{e}\v{z}\'{\i}~4, CZ-182 08 Prague 8, Czechia,}
\affil[3]{Faculty of Civil Engineering, The Czech Technical University, Th\'akurova 7, CZ-166 29 Prague, Czechia,}
\affil[4]{Institute of Mathematics, Faculty of Science, University of South Bohemia, Brani\v sovsk\' a~31, CZ-37005 \v{C}esk\'{e} Bud\v{e}jovice,  Czechia}

\begin{document}

\maketitle

\begin{abstract}
A sharp-interface model describing static equilibrium configurations of shape memory alloys by means of interfacial polyconvex energy density introduced in \cite{Silhavy-2010} and extended to a quasistatic situation in \cite{hkmk} is  computationally tested.  Elastic properties of  variants of martensite and the austenite are described by polyconvex energy density functions. Volume fractions of particular variants are  modeled by a map of bounded variation. Additionally, energy stored in martensite-martensite and austenite-martensite interfaces is measured by an interface-polyconvex function. It is assumed that transformations between  material variants are  accompanied by energy dissipation which, in our case, is positively and one-homogeneous giving rise to a rate-independent model. Various two-dimensional computational examples are  presented and the used computer  code is made available for downloads.
 
\end{abstract}

\section{Introduction}
Reversible solid-to-solid phase transformations in shape memory alloys (abbreviated ``SMA'' in the sequel) give rise to well-known shape memory effects, which have various technological applications \cite{Mohd-Jani-2014}. Such materials have a high-temperature phase called austenite and a low-temperature phase called martensite. The austenitic phase has only one variant but the martensitic phase exists in many symmetry related variants and can form a microstructure by mixing those variants (possibly also with austenite) on a fine scale \cite{Bhattacharya-2003}. The most common examples of SMA are, e.g., Ni-Ti, Cu-Al-Ni or In-Th.

In the last three decades, many models of SMA varying in modeling approach, scale and purpose appeared in the literature; see \cite{Ben-Zineb-2016} for a recent reviews. Single crystal models often attempt a plausible description of formation and/or evolution of fine microstructures of various types \cite{Kruzik-2005,Stupkiewicz-2007,Seiner-2011,Tuma-2016a}.

Variational models for SMA microstructures assume that the formed structure has some optimality property. The reason for the formation of microstructures is that no exact optimum can be achieved and optimizing sequences have to develop finer and finer oscillations. The goal is to model presence of different phases, which leads to the so-called  multi-well structure of material stored energy density. If the temperature $\vartheta$ is below the transformation temperature $\vartheta_t$  then the stored energy density  is minimized on wells SO$(3)U_i$, $i=1,\ldots, M$,  defined by  $M$ positive definite and symmetric  matrices $U_1,\ldots, U_M$. Above the transformation temperature, the global minimizer of the energy density is just the special orthogonal group SO$(3)$, describing the stress-free strain of austenite.  Naturally, at the transformation temperature, all martensitic variants as well as the austenite must be considered and the stored energy density is minimized on  $M+1$ wells where the well of the austenite is just SO$(3)$. This behavior generically causes nonexistence of minimizers even in elastostatic problems.  A way out is relaxation in the calculus of variations searching for the so-called quasiconvex envelope of the specific stored energy \cite{dacorogna-book, mueller-1999} or using Young measures \cite{k-p1,k-p,mielke-roubicek-2003}. The downside of these techniques is that we do not have a closed formula of the envelope at our disposal and that  physically justified conditions on deformations as orientation-preservation and injectivity are not included in these models. On the other hand, recently a few new results appeared and we cite  \cite{benesova-kruzik-2017} for a survey. We also refer to \cite{benesova-kruzik-2015} for a weak* lower semicontinuity results for sequences of bi-Lipschitz orientation-preserving maps in the plane and to \cite{benesova-kampschulte-2014} for an analogous result along  sequences of quasiconformal maps. Then \cite{koumatos-rindler-wiedemann-2013} found relaxation including orientation preservation for $1<p\le n$, where $n$ is a dimension of the problem  and $p$ is the power with which the deformation gradient is integrable. Finally, \cite{Conti-Dolzmann-2015} derived a relaxation result for orientation preserving deformations with an extra assumption on the resulting functional. Let us point out that the right Cauchy-Green strain tensor $F^\top F$ maps SO$(3)F$ as well as (O$(3)\setminus$SO$(3)$)$F$ to the same point. Here O$(3)$ are orthogonal matrices with  determinant $\pm 1$. Thus, for example, $F\mapsto |F^\top F-\mathbb{I}|$, where $\mathbb{I}$ is the identity matrix,  is minimized on two  energy wells, i.e., on SO$(3)$ and also on O$(3)\setminus$SO$(3)$.  However, the latter set is not acceptable in elasticity. Additionally, notice that, for example, considering arbitrary $Q\in {\rm O}(3)\setminus {\rm SO}(3)$ and an arbitrary $R\in {\rm SO}(3)$ such that $Q$ and $R$ are rotations around the same axis of the Cartesian system then rank$(Q-R)=1$.  Consequently, corresponding effective macroscopic energy density (relaxation or quasiconvexification) resulting from the original multiwell energy density can be smaller then a physically relevant energy density.  

A well justified model of materials in nonlinear elasticity which allows for orientation preservation and injectivity is based on polyconvexity due to J.M.~Ball \cite{ball-1977,ball-1981}. It is also relatively easy to construct polyconvex energy density. On the other hand, multi-well stored energy densities generically used to model shape memory alloys are not polyconvex and accompanied with  non-existence of  minimizers of the total energy in Sobolev spaces. Recently a new static model of such materials appeared in \cite{Silhavy-2010,Silhavy-2011} which combines polyconvex energy densities  of austenite and martensitic variants with a newly defined  polyconvexity of interfacial energy between variants of martensite and austenite. Volume fractions of martensitic variants and of austenite play a role of  additional design variables. Let us note here that it is well agreed that the energy (alternatively called interfacial, surface, interaction, etc.) stemming from the presence of interfaces in the material is crucial for determining the proper size scale in the microstructures \cite{Niclaeys-2002,Govindjee-2007,Petryk-2010,Bernardini-2014,Tuma-2016b}. However, as demonstrated by Seiner et al. \cite{Seiner-2011}, the real micromorphology does not necessary correspond to the minimum of the mere sum of the elastic (bulk) and interface (surface) energies obtained in static models.

In this work, we consider the model from \cite{hkmk} which extended  the static model of \v{S}ilhav\'{y} to a rate-independent evolutionary model and proved existence of an energetic solution. Obtained time dependent deformations are orientation preserving and injective and no additional regularization of variables is needed if passing from a static to an evolutionary model. Moreover, each martensitic variant and the austenite are assumed to be polyconvex and obey realistic properties, see \eqref{W-ass1}--\eqref{W-ass4} below. The main aim of this work is to provide computational examples with the rate-independent model. Moreover, a Matlab code  used for calculating our examples is freely provided at \url{https://www.mathworks.com/matlabcentral/fileexchange/68547}.

The plan of the paper is as follows. After introducing used notations, we  first  describe our model, the stored energies, dissipation and loading. Then we recall the existence of an energetic solution proved in \cite{hkmk}. As it is nowadays a standard procedure, cf.~e.g., \cite{FrancfortMielke-2006,MielkeRoubicek-2015,MielkeTheil-1999}, we only sketch main steps. Finally, we demonstrate viability of the complete modeling concept on   numerical simulations.

\medskip

\section{Notation and Preliminaries}
 Let $d=2$ or $d=3$. By $\mathcal{L}^d$, we denote the $d$-dimensional Lebesgue measure and by $\HH^{d-1}$ the $d-1$-dimensional Hausdorff measure.  For two matrices $A = (a_{ij}), B = (b_{ij}) \in \R^{d\times d}$, we define $A:B = a_{ij}b_{ij}$ with Einstein's sum
convention.  By $(A\times n) \in\R^{d\times d}$ we denote the tensor defined by $(A \times n) b = A (n \times b)$, i.e. $(A \times n)_{kj} = \eps_{\ell ij}
a_{k\ell} n_i$, where $\eps_{\ell ij}$ is the Levi-Civita symbol.

\medskip

For $A \in \R^{d\times d}$ invertible , we have
$$\cof A=(\det A) A^{-\top}\ .$$
In particular,  if $d=2$ then 
\begin{align}\label{cof2d}
\cof F:= 
\begin{pmatrix} 
F_{22} & -F_{21} \\
-F_{12} & F_{11} 
\end{pmatrix}\ ,
\end{align}
where $F_{ij}$ are entries of $F$.

We refer e.g. to \cite{Gurtin-Struthers-1990} for a definition of the surface gradients $\nabla _S$. If $n\in\R^d$ is an outer unit normal to the surface $S$, then 
\begin{align} \label{surface_deformation_gradient}
\nabla_S:=\nabla(\mathbb{I}-n\otimes n), 
\end{align}
where $\mathbb{I}$ denotes the unit matrix in $\R^{d\times d}$. We denote by $W^{1,p}$ a standard Sobolev space of maps which are together with their distributional derivatives integrable with the $p$-th power, cf., e.g.,  \cite{EvansGariepy-Book}. Further, BV stands for the space of integrable maps with bounded variation. We refer, e.g., to \cite{Ambrosio-Fusco-Pallara-Book} for a detailed exposition on this subject.  
The space of vector-valued Radon measures on $\O$ with values in $\R^{N}$  will be denoted ${\cal M}(\O;\R^{N})$.

Let $\tilde\O\subset\O\subset\R^d$, $d>1$, be Lebesgue measurable sets with a finite Lebesgue measure and let $B(x,r):=\{a\in\R^d:\, |x-a|< r\}$. If $x\in\O$ we define 
$$\theta(\tilde\O,x):=\lim_{r\to 0}\mathcal{L}^d(\tilde\O\cap B(x,r))/\mathcal{L}^d(B(x,r))$$ whenever this limit exists. We call it  the density of $\tilde\O$ at $x$. We call  $\{x\in\O;\, \theta(\tilde\O,x)=1\}$ the set of points of density of $\tilde\O$.
If $\theta(\tilde\O,x)=0$ for some $x\in\O$ we call $x$ the point of rarefaction of $\tilde\O$. The measure-theoretic boundary  of $\tilde\O$ denoted $\partial^*\tilde\O$ is the set of all points $x\in\O$ such that $\theta(\tilde\O,x)\ne 0$, nor 1,  or  $\theta(\tilde\O,x)$ does not exist. We call $\tilde\O$ a set of finite perimeter if $\mathcal{H}^{d-1}(\partial^*\tilde\O)<+\infty$. Let $n\in\R^d$ be a unit vector and let $H(x,n):=\{\tilde x\in\O:\, (\tilde x-x)\cdot n<0\}$. We say that $n$ is a measure-theoretic normal to $\tilde\O $ at $x$ if $\theta(\tilde\O\cap H(x,-n),x)=0$ and $\theta((\O\setminus\tilde\O)\cap H(x,n),x)=0$.  The measure-theoretic normal exists for $\HH^{d-1}$ almost every point in $\partial^*\tilde\O$. More details can be found in \cite{EvansGariepy-Book, Silhavy-1997}. 

\section{Model description}

\subsection{Elastic energy} \label{ss-static} %

{\it Geometry and total energy: } We assume that the specimen in its reference configuration is represented by a bounded
Lipschitz domain $\Ome\subset\R^d$. We consider a shape memory alloy which allows for $M$ different variants of martensite. We assume that the region occupied by the $i$-th variant of martensite is given by $\Ome_i \subset \Ome$ for
$1 \leq i \leq M$, while the region occupied by austenite is given by $\Ome_0 \subset \Ome$. In particular, the sets
$\Ome_i$ are pairwise disjoint for $0\le i\le M$. We also assume that the sets $\Ome_i$ are open, that $N := \Ome \BS \bigcup_i \Ome_i$
is a zero (Lebesgue) measure set and that each $\Ome_i$ is of a finite perimeter, i.e., that its characteristic function has a bounded variation. The case $\Ome_i=\emptyset$ for some $0\leq i\leq M$ is not excluded. The partition of $\Ome$ into
$\{\Ome_i\}_{i=0}^M$ can be then identified with a mapping $z:\O\to\R^{M+1}$ such that $z_i(x)=1$ if $x\in\Ome_i$ and
$z_i(x)=0$ else. Clearly, $\sum_{i=0}^Mz_i(x)=1$ for almost every $x\in\Ome$.  We will call $z$ the partition map corresponding to $\{\Ome_i\}_{i=0}^M$.  
 We hence consider $z \in \ZZ$, where
\begin{multline}
  \ZZ := \Big \{ z \in {\rm BV}(\Ome,\R^{M+1}) : z_i\in\{0,1\}\, ,  z_i z_j = 0 \text{ for $i \neq j$, \ }\nonumber \\ \sum_{i=0}^M z_i(x) = 1 \text{ for 
    a.e $x \in \Ome$} \Big \}.
\end{multline}
In order to describe the state of the elastic material completely, we also need to introduce the deformation function $v
\in W^{1,p}(\Ome,\R^d)$, $p > d$, which describes the deformation of the elastic body with respect to the reference
configuration $\Ome$. We assume that part of the boundary $\Gam_{\rm Dir} \subset \p \Ome$ with $\HH^{d-1}(\Gam_{\rm Dir}) \neq 0$ is
fixed. We hence assume that there is $y_{\rm Dir}\in W^{1,p}(\O;\R^d)$ such that $y=y_{\rm Dir}$ on $\Gamma_{\rm Dir}\subset\partial\Ome$, i.e. $y_{\rm Dir}$ represents a Dirichlet boundary condition. We hence consider deformations $v \in \VV$, where
\begin{align}
  \VV  = \Big\{ v \in W^{1,p}(\Ome,\R^d) : \det \nabla v>0\ ,
	 \int_\O\det\nabla v(x)\,\md x\le \mathcal{L}^d(v(\O)) \Big\},
\end{align}
where we will always use the assumption $p > d$. The integral inequality together with the orientation-preservation is the so-called Ciarlet-Ne\v{c}as condition which ensures invertibility of $v$ almost everywhere in $\O$ \cite{Ciarlet-Book,CiarletNecas-1987}. In the following, we will assign to each state of the material $(v,z)
\in \VV \times \ZZ$ an energy $E(v,z)$. In our model, the energy consists both of a bulk part, penalizing deformation
within the single phases, and an interfacial energy, penalizing deformation of the interfaces between the phases.

\medskip

{\it Bulk energy:} The total bulk energy of the specimen has the form
\begin{align} \label{stat-bulk}  %
  E_{\rm b}(y,z)=\int_\Ome W(z(x),\nabla y(x))\ \md x.
\end{align}
We assume that the specific energy $W$ of the specimen can be written as
\begin{align}
  W(z(x),F) := z(x)\cdot \hat W(F) := \sum_{i=0}^M z_i(x) \hat W_i(F),
\end{align}
where $\hat W_i$, $0 \leq i \leq M$ is the specific energy related to the $i$-th phase of the material and $\hat W:=(\hat W_0,\ldots, ,\hat W_M)$. We will work in the
framework of hyperelasticity, where the first Piola-Kirchoff stress tensors of austenite and martensite have polyconvex
potentials denoted by $\hat W_0$ (austenite) and $\hat W_i$, $i=1,\ldots, M$ for each variant of martensite.  For $0\le
i\leq M$, we therefore assume
\begin{align} \label{W-ass1} %
  \hat W_i(F):=
  \begin{cases}
    h_i(F,\cof F,\det F) & \mbox{ if } \det F>0, \\
    +\infty \mbox{ otherwise.}
  \end{cases}
\end{align}
for convex functions $h_i:\R^{19}\to\R$ if $d=3$. and $h_i:\R^9\to\R$ if $d=2$. If $d=2$ then it follows from \eqref{cof2d} that $\cof F$ is redundant in the definition of polyconvexity  but we keep it there for simplicity. 
We  use the following additional standard assumptions on the specific bulk energies $\hat W_i$. For $0\leq i\leq M$ and
$F\in\R^{d\times d}$, we assume that for some $C>0$ and $p>d$ 
\begin{align} 
  &\hat W_i(F)\ge C(-1+|F|^p) &&\forall F \in \R^{d \times d}\ , \label{W-ass2} \\
  &\hat W_i(RF)=\hat W_i(F) && \forall R\in\SO(d), F \in \R^{d\times d} \ , \label{W-ass3}\\
	&\lim_{\det F\to 0_+} \hat W_i(F)=+\infty\ . \label{W-ass4}
\end{align}

\medskip

{\it Interfacial energy:} We will consider the interfacial energy as the one introduced by \v{S}ilhav\'{y} in
\cite{Silhavy-2010,Silhavy-2011}. We therefore assume that the specific interfacial energy $f_{ij}$  between the two different phases $i, j\in\{0,\ldots, M\}$
can be written in the form
\begin{align} \label{intcon-1} %
  \frac12 f_{ij}(F, n)=g_i(F,n)+g_j(F,n),
\end{align}
where $F\in\R^{d\times d}$ and $n\in\R^d$ is a unit vector such that $F n=0$.
We assume
\begin{align} \label{intcon-2} %
  g_i(F,n):= \begin{cases}\Psi_i(n,F\times n, \cof F\, n) & \text{ if $d=3$,}\\
  \Psi_i(n,\cof F\, n) & \text{ if $d=2$.}
  \end{cases}
\end{align}
where the functions $\Psi_i:\R^{15}\to\R$ if $d=3$ and  $\Psi_i:\R^{4}\to\R$ are nonnegative convex and positively one-homogeneous for $i=1,\ldots,
M$.  As in \cite{Silhavy-2010}, we assume for $0 \leq i \leq M$, $ \forall F\in\R^{d\times d}$, $\forall n\in S^2$ 
\begin{align}
  g_i(RF,n)& =g_i(F,n)\  \forall R\in \SO(d)\ ,  \label{g-ass2} \\
  g_i(F,n) & =g_i(F,-n)\ .  \label{g-ass3}
\end{align}
 We introduce a subspace $\QQ \subset \VV \times \ZZ$ of functions with ``finite
interfacial energy'', using a slightly modified version of \cite[Def. 3.1]{Silhavy-2010}. It is given as follows:

\begin{definition}[states with finite interfacial energy] \label{def-intok} %
 Let $d=3$.  For any pair $y,z \in \VV \times \ZZ$ let $S_i=\partial^*\Ome_i\cap\Ome$ where $\partial^*\Ome_i$ is the
  measure-theoretic boundary of $\Ome_i$ with outer (measure-theoretic) normal $n_i$. We denote by $\QQ \subset \VV
  \times \ZZ$ the set of all pairs $(y,z) \in \VV \times \ZZ$ such that for every $0 \leq i \leq M$ there exists a
  measure $J_i:=(b_i,H_i,c_i)\in{\cal M}(\O;\R^{15})$ with $ b_i := n_i \HH^{d-1}_{|S_i}$ and 
  \begin{align} \NT %
    \int_{\Ome_i} \nabla y\ (\nabla \times v) \ \md x=\int_\Ome v\ \md H_i \ %
    &&\text{and}&& %
    \int_{\Ome_i} (\cof\nabla y) : \nabla v \ \md x=\int_\Ome v\cdot \md c_i
  \end{align}
  for all $v\in C^\infty_0(\O;\R^d)$. If $d=2$ we define $J_i:=(b_i,c_i)\in {\cal M}(\O;\R^{4})$.
\end{definition}
\begin{remark}\label{measures}
\upshape
  The measures $H_i$ and $c_i$ can also be expressed for smooth $y$  as
  \begin{align}\label{measures1}
     b_i := n_i \HH^{d-1}_{|S_i}, && %
     H_i:= \nabla_{S_i} y\times n_i \HH^{d-1}_{|S_i} && \text{and } && %
     c_i:= (\cof\nabla_{S_i} y) n_{|S_i}
  \end{align}
  for all $v\in C^\infty_0(\O;\R^d)$. Indeed, in components
  and using Einstein's summing convention, we have
  \begin{align*}
    \int_\Ome [\nabla y (\nabla \times v)]_k %
    &=   \int_\Ome [\p_j y_k \eps_{j\ell m}  \p_\ell v_m] dx \\ 
    &=   - \int_\Ome [\p_\ell \p_j y_k \eps_{j\ell m}  v_m]  dx + \int_{\p \Ome} [n_\ell \p_j y_k \eps_{j\ell m}  v_m ]  dx \\
    &=   \int_{\p \Ome} [\nabla v]_{kj} [n \times v]_j  dx  =  \int_{\p \Ome} [(\nabla y \times n) v]_k dx,
  \end{align*}
  where we used the  identity $\nabla \times \nabla y = 0$.  With the notation $(\cof \nabla y)_{ij} = b_{ij}$, we also have
  \begin{align*}
    \int_\Ome (\cof \nabla y) : (\nabla v) dx %
    &=   \int_\Ome [b_{kj} \p_j v_k] dx \\
    &=   - \int_\Ome [\p_j b_{kj}   v_k]  dx + \int_{\p \Ome} [n_j b_{kj}  v_k ]  dx \\
    &= \int_{\p \Ome} [(\cof \nabla y) n]_{k} v_k dx = \int_{\p \Ome} (\cof \nabla y) n \cdot v \ \md x.
  \end{align*}
  where we used the Piola identity $\nabla \cdot (\cof \nabla y) = 0$.  
\end{remark}
With the notation of Definition \ref{def-intok}, we define the interfacial energy as
 \begin{align} \label{stat-interface}
   E_{\rm int}(y,z):=
   \begin{TC}
     \displaystyle \sum_{i=0}^M \int_\Ome \Psi_i\left(\frac{\md J_i}{\md |J_i|}\right)\ \md|J_i|\ %
     &\text{for $(y,z) \in \QQ$,} \\
     \infty &\text{else.}
   \end{TC}
 \end{align}
Here $|J_i|$ denotes the total variation of the measure $J_i$. 
We assume that for all $0\le i\le M$ there is some $c>0$  that for all $A$
\begin{align}\Psi_i(A)\ge c|A| . \label{g-ass1} \end{align}

\begin{remark}
Notice that in view of \eqref{measures1} and the definition of $b_i$ \eqref{g-ass1} allows us to bound the measure-valued derivative of $z$,  $\|D z\|_{\mathcal{M}(\O;\R^{(M+1)\times 3})}$,  by means of    $E_{\rm int}(y,z)$.  Consequently, we control $\|z\|_{{\rm BV}(\O;\R^{M+1})}$ in terms of the interfacial energy.
\end{remark}

\medskip

Examples of  surface energy densities include for example \cite{Silhavy-2011} $g_i(F,n)=\alpha|F|=\alpha|F\times n|$ for $\alpha>0$ or $g_i(F,n)=\alpha |\cof F n|$.

{\it Body and surface loads:} We assume that the body is exposed to possible body and
surface loads, and that it is elastically supported on a part $\Gamma_0$ of its boundary. The part of the energy related to this loading is given by a 
 functional $L\in C^1([0,T]; W^{1,p}(\O;\R^d))$.   We  consider 
\begin{align} \label{staticloading} %
  L(t,y) &:=-\int_\Ome b(t,x)\cdot y(x)\ \md x-\int_{\Gamma_1} s(t, x)\cdot y(x)\ \md A\nonumber\\ &+\frac{K}{2}\int_{\Gamma_0}|y(x)-y_D(t,x)|^2\,\md A\ ,
\end{align}
where $b(t,\cdot):\O\to\R^d$ represents volume density of body forces and $s(t,\cdot):\Gamma_1\subset\partial\O\to\R^d$ is areal density of surface forces 
applied on a part $\Gamma_1$ of the boundary. The last term in \eqref{staticloading}  with $y_D(t,\cdot)\in W^{1,p}(\O;\R^d)$  represents  energy of a spring with a spring stiffness  constant $K>0$. Thus, our specimen is elastically supported on $\Gamma_0$ in such a way, that for $K\to\infty$ 
$y$ is forced to be close to $y_D$ on $\Gamma_0$ in the sense of the $L^2(\Gamma_0;\R^d)$ norm.  A term of this type already appeared in \cite{kruzik-otto-2004}  and its static version also in  \cite{mielke-roubicek-2003}.  Namely, prescribing a boundary condition  from $W^{1-1/p,p}(\partial\Ome;\R^d)$ \cite{marschall},  it is generally not known whether it can be extended to the whole $\Omega$ in such a way  that the extension lives in $\VV$. It is, to our best knowledge,  an unsolved problem in three dimensions and therefore it is generically assumed in nonlinear elasticity  that such an extension exists; cf.~\cite{Ciarlet-Book}, for instance. 
The last term in \eqref{staticloading} overcomes this drawback. Namely, if $y_D$ cannot be extended  from the boundary as an orientation-preserving map the term in question can never be zero.   
\medskip


\subsection{Dissipation}
Since the evolution in SMA is typically connected with energy dissipation, we need to define a suitable dissipation function. Experimental evidence shows that considering a rate-independent dissipation mechanism is a reasonable approximation in a wide range of rates of external loads, hence, this dissipation is to be positively one-homogeneous. We simply associate the dissipation to the magnitude of the time derivative of
$z$, i.e., to $|\dot z|_{M+1}$, where $|\cdot|_{M+1}$ is a norm on the $\R^{M+1}$ space.\footnote{For more refined rate-independent dissipation functions tailored for SMA polycrystals see, e.g. \cite{Junker-2014,Frost-2016}.} Therefore, the specific dissipated energy associated with a change of the variant distribution from $z^1$ to $z^2$ is postulated as
\begin{align}\label{dissipation}
  D(z^1,z^2):=|z^1-z^2|_{M+1}.
\end{align}
Then the total dissipation reads
\begin{align*}
  \mathcal{D}(z^1,z^2):=\int_\Ome D(z^1(x),z^2(x))\ dx \ .
\end{align*}
The $\DD$-dissipation of a curve $z(t)$ is defined as
\begin{align*}
  {\rm Diss}_\DD(z,[s,t]) = \sup \Big \{ \sum_{j=1}^N \DD(z(t_{i-1}),z(t_i)) : N \in \N, s = t_0 \leq \ldots \leq t_N = t \Big \}
\end{align*}

We denote for $(t,y,z)\in[0,T]\times\VV\times\ZZ$ the total energy 

\begin{align}
\EE(t,y,z):=E_{\rm b}(y,z)+E_{\rm int}(y,z)-L(t,y)\ .
\end{align}

\subsection{Energetic solution} 

\medskip

Suppose, that we look for the time evolution of $t\mapsto y(t)\in \VV$ and $t\mapsto z(t)\in \ZZ$ during a process time interval $[0,T]$ where  $T>0$ is the time horizon. We use the following notion of solution from \cite{FrancfortMielke-2006}, see also
\cite{MielkeTheil-2004,MielkeTheilLevitas-2002}: For every
admissible configuration, we ask the following conditions to be satisfied for all $t \in [0,T]$.
\begin{definition}[Energetic solution]
  We say that $(y,z) \in \VV\times \ZZ$ is an energetic solution  to $(\EE,\DD,L)$ if $t \mapsto \p_t E(y(t),z(t)) \in L^1((0,T))$
  and if for all $t \in [0,T]$, the stability condition (S) and the condition (E) of energy balance are satisfied,
  where
\begin{align}
  &\EE(t, y(t),z(t))\leq \EE(t,\tilde y,\tilde z)+\mathcal{D}(z(t),\tilde z) \qqquad %
  \text{$\forall (\tilde y,\tilde z)\in\QQ$}. \tag{S}
\end{align}
and 
  \begin{align}  \tag{E}      %
    \begin{aligned}
    &\EE(t, y(t),z(t))+{\rm Diss}_\DD(z;[0,t]) \\
    &\qqquad=\EE(0,y(0),z(0)) +\textstyle \int_0^t \frac{\partial \EE}{\partial
      t}(s, y(s),z(s))\, \md s \hspace{-2ex}   
    \end{aligned} 
  \end{align}
  are satisfied.
\end{definition}

An important role is played by the so-called stable states defined for each $t\in[0;T]$. We set
$$
\mathbb{S}(t):=\{(y,z)\in\VV\times \ZZ:\, \EE(t, y,z)\leq \EE(t,\tilde y,\tilde z)+\mathcal{D}(z,\tilde z)\,\forall (\tilde y,\tilde z)\in\QQ\}\ . $$

\section{Existence of the energetic solution}

A standard way how to prove the existence of an energetic solution is to construct time-discrete minimization problems
and then to pass to the limit. Before we give the existence proof, we need some auxiliary results. If $N\in\N$ define time increments  
$t_k:=kT/N$ for $0\le k\le N$. Further, we abbreviate  $q:=(y,z)\in\QQ$.  Assume that at $t=0$  there is given an initial distribution of 
phases $z^0\in\ZZ$ and $y^0\in\VV$  such that $q^0=(y^0, z^0)\in \mathbb{S}(0)$.   For $k=1,\ldots, N$,  we define a sequence of minimization problems 
\begin{align}\label{incremental}
\text{minimize } \EE(t_k, y,z)+\DD(z,z^{k-1})\ ,\ (y,z)\in \QQ\ .
\end{align}

Denoting $B([0,T];\VV)$ the set of bounded maps $t\mapsto y(t)\in\VV$ for all $t\in[0,T]$, we have the following result showing the existence of an energetic solution proved in \cite{hkmk}.

\begin{theorem}
Let  $T>0$, $p>d$, $y_{\rm D}\in C^1([0,T]; W^{1,p}(\O;\R^d))$, \eqref{W-ass1}-\eqref{W-ass4}, \eqref{intcon-2}, \eqref{g-ass3}-\eqref{g-ass1}. Let  $(y(0),z(0))\in S(0)$ and Let there be $(y,z)\in\QQ$ such that $\EE(0,y,z)<+\infty$.  Then there is and energetic solution to the problem $(\EE,\DD,L)$ such that $y\in B([0,T];\VV)$, $z\in {\rm BV}([0,T]; L^1(\O;\R^{M+1}))\cap L^\infty(0,T;\ZZ)$. 
\end{theorem}

\section{Computational illustrations}
\begin{figure}
\centering 
 \includegraphics[width=\textwidth]{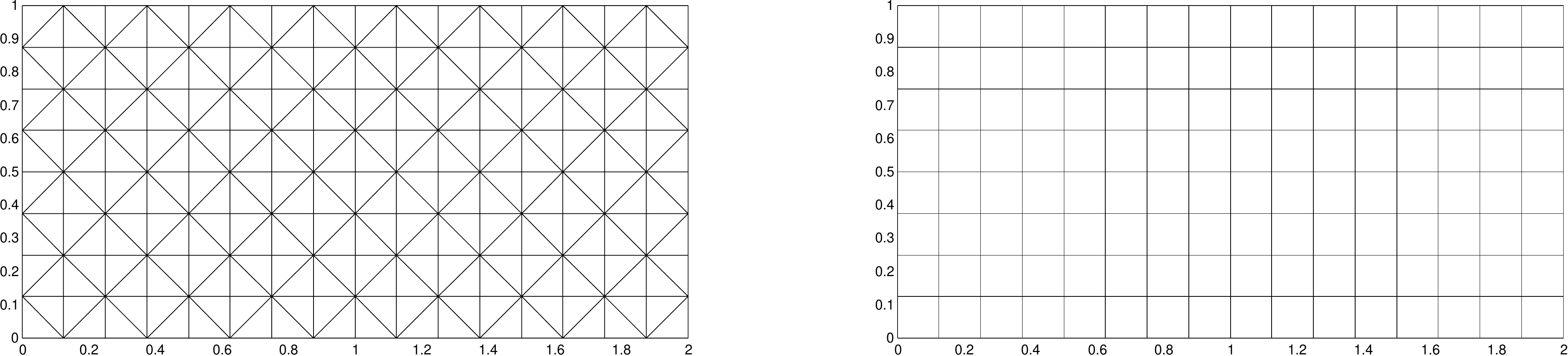} 
\caption{A FEM triangular mesh $\mathcal{T}$ (left) and the underlying rectangular mesh  (right)  for visualization.}
\label{mesh}
\end{figure}

\begin{figure}[t]
\centering 
 \includegraphics[width=0.5\textwidth]{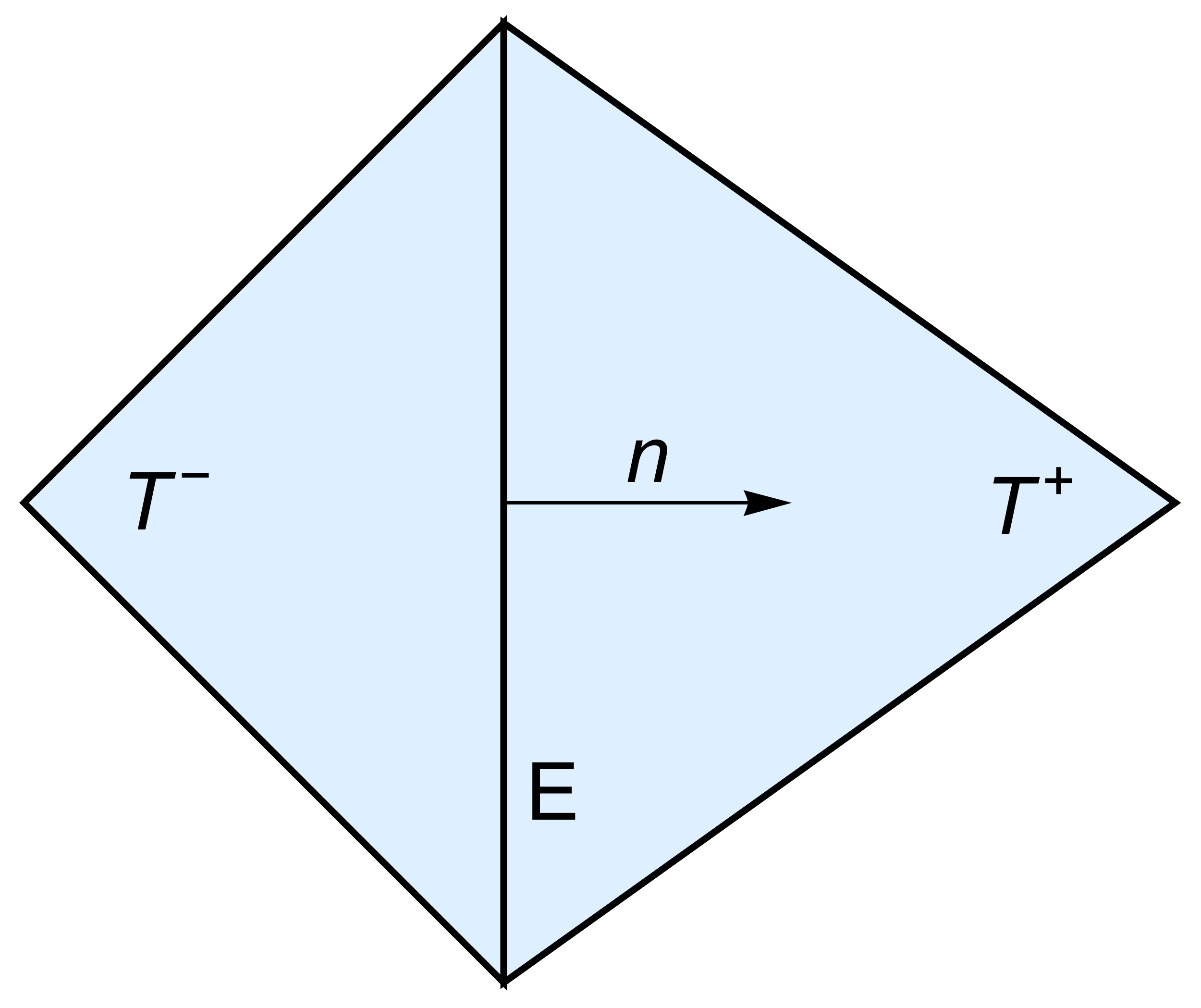} 
\caption{An internal edge $E \in \mathcal{E_I}$ shared by two elements $T^+, T^-  \in \mathcal{T}$ and its normal vector $n$. The surface and interfacial energies are evaluated on $E$, while the bulk energy and the dissipation on $T^+$ and $T^-$.}
\label{edge}
\end{figure}

In the following two-dimensional $(d=2)$ illustrative examples, we deal with SMA subjected to loading at temperature $\vartheta < \vartheta_{t}$, i.e. martensite is the stable phase. We consider two martensitic variants, hence $z$ can be condensed to a simple scalar variable. Moreover, we assume that $K=+\infty$ in \eqref{staticloading}, i.e., Dirichlet boundary conditions, and   we also  assume that the  specimen occupies a rectangular domain $$\Omega= (0,2) \times (0,1)\ .$$ A regular triangulation $\mathcal{T}$ of $\overline \Omega$ (an example is given in Figure \ref{mesh}) is considered with $\mathcal{E}$ denoting the set of all edges in $\mathcal{T}$ and $\mathcal{E_I}$ its subset of internal edges. Given an edge $E \in  \mathcal{E_I}$, we can assign two triangles $T^+ \in \mathcal{T}$ and $T^-  \in \mathcal{T}$ such that $E \subseteq \partial T^+   \cap  \partial T^- $ (eg. the edge $E$ is shared by the triangles  $T_+$ and $T_-$). A vector $n \in \R^2$ such that $ |n|=1$ denotes a normal vector orthogonal to $E$, see Figure \ref{edge}. The symbol $\mathcal{N}$ denotes the set of all nodes. \\

\begin{figure}
\centering 
 \includegraphics[width=\textwidth]{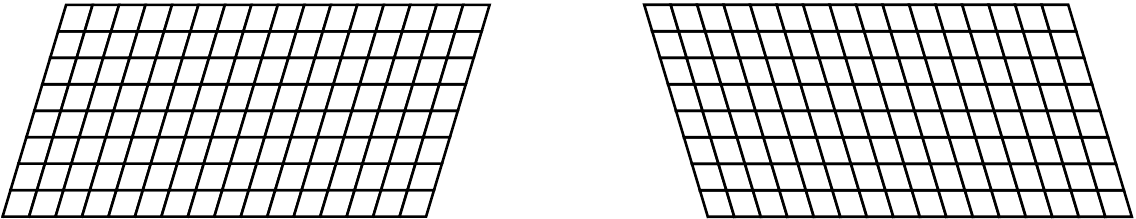} 
\caption{Examples of mesh deformations corresponding to stretching matrices $F_1$ (left), $F_2$ (right) for $\epsilon=0.3$. }
\label{mesh_stretching}
\end{figure}

We consider the interfacial energy $E_{\rm int}$ in the form
$E_{\rm int}= E^1_{\rm int}+E^2_{\rm int}$. Namely, 
\begin{align}
 E^1_{\rm int}(y,z):= \int_{E \in  \mathcal{E_I}} \alpha_{\rm i} | z^+ (s) -   z^-(s) | \ \md s,
\end{align}
where $ z^+$ and $ z^-$ denote restrictions of $z$ to triangles $ T^+$ and $ T^-$ and $\alpha_{\rm i} \geq 0$ is a parameter, and 
\begin{align}E^2_{\rm int}(y,z):= \int_{E \in  \mathcal{E_I}}  
\alpha_{\rm s}
|\cof \mathbb{F}(s) n | | z^+ (s) -   z^-(s) | \ \md s
\end{align}
is based on the cofactor of the surface deformation gradient $\mathbb{F}$ given by \eqref{surface_deformation_gradient}.
A parameter $\alpha_{\rm s} \geq 0$ is given. The bulk energy $E_{\rm b}$ is considered in the form \eqref{stat-bulk} with the specific energy
\begin{align}
  W(z(x),F) =  z(x) \hat W_1(F) + (1-z(x)) \hat W_2(F).
\end{align}
Here, $\hat W_1, \hat W_2$ are densities in the form
\begin{equation*}
\hat W_1(F):=  \underline{W}(F F_1^{-1}), \quad \hat W_2(F):= \underline{W}(F F_2^{-1}),
\end{equation*}
where $F_1, F_2$ are given stretching matrices
\begin{align}
F_1:= 
\begin{pmatrix} 
1 & \epsilon \\
0 & 1 
\end{pmatrix}, \quad 
F_2:= 
\begin{pmatrix} 
1 & -\epsilon \\
0 & 1 
\end{pmatrix}
\end{align}
defined by a parameter $\epsilon > 0$. The form of  \underline{W} is given by the compressible two-dimensional Mooney-Rivlin material model
\begin{eqnarray}\label{density_Mooney-Rivlin_F}
\underline{W}(F):=\alpha \tr (F^T F)  + \delta_1 (\det F)^2 - \delta_2 \ln (\det F).
\end{eqnarray}
Parameters $\alpha, \delta_1, \delta_2 > 0$ satisfy the relation\footnote{In terms of invariants, i.e. eigenvalues of $F^\top F$ we have that $$\hat W(\lambda_1,\lambda_2)=\alpha (\lambda_1^2+\lambda_2^2)+\delta_1\lambda_1^2\lambda_2^2-\delta_2\ln\lambda_1\lambda_2\ .$$ Putting the first derivatives to  zero yields
\begin{eqnarray*}
\lambda_1\partial \hat W/\partial\lambda_1= 2\alpha\lambda_1^2+2\delta_1\lambda_1^2\lambda_2^2-\delta_2=0, \\
\lambda_2\partial \hat W/\partial\lambda_2= 2\alpha\lambda_2^2+2\delta_1\lambda_1^2\lambda_2^2-\delta_2=0.
\end{eqnarray*}
Subtracting  the two equations yields $\lambda_1=\lambda_2$ (both must be $>0$).
If $2\alpha+2\delta_1=\delta_2$ we get that $\lambda_1=\lambda_2=1$ minimizes $\hat W$.
} 
\begin{equation}
\delta_2 = 2 \alpha + 2 \delta_1
\end{equation}
and it holds
$ \underline{W}(F) \rightarrow \infty \mbox{ for } \det F \rightarrow 0+$. Since $I = \mbox{argmin}_{F} \, \underline{W}(F)$, it holds
$$ F_1 = \mbox{argmin}_{F} \, \hat W_1(F), \quad F_2 = \mbox{argmin}_{F} \, \hat W_2(F). $$ Therefore, the minimizer of the bulk energy $E_{\rm b}$ aligns the deformation $F$ with one of the stretching matrices: $F_1$ (for $z=0$ everywhere in $\Ome$) or $F_2$ (for $z=1$ everywhere in $\Ome$), see Figure
\ref{mesh_stretching}.

\begin{remark}
It is convenient to reformulate \eqref{density_Mooney-Rivlin_F} in terms of Green's deformation tensor $C=F^T F$ as 
\begin{eqnarray}\label{density_Mooney-Rivlin_C}
\underline{\underline{W}}(C):=\alpha \tr (C)  + \delta_1 \det C - \frac{\delta_2}{2} \ln (\det C)
\end{eqnarray} 
 and exploit relations
\begin{eqnarray}\label{transformation_isotropic_C}
\hat W_1(F) = \underline{\underline{W}}(F_1^{-T} C  F_1^{-1}), \quad 
\hat W_2(F) = \underline{\underline{W}}(F_2^{-T} C  F_2^{-1}).
\end{eqnarray}
For anisotropic energy densities, we should replace $F_1^{-T} C  F_1^{-1}, F_2^{-T} C  F_2^{-1}$ in \eqref{transformation_isotropic_C} by 
$$U_1^{-T} C  U_1^{-1}, \quad U_2^{-T} C  U_2^{-1},$$ where $U_1^2=F_1^T F_1$ and $U_2^2=F_2^T F_2$. 
\end{remark}
The dissipation is assumed in the form 
\begin{align}
\DD(z,z^{k-1}) := \int_\Ome \beta | z - z_{k-1} | \ \md x, 
\end{align}
where $\beta \geq 0$ is a given parameter.  \\

The finite element method (FEM) is applied for the discretization of all functionals above. We choose the lowest possible order finite element functions: 
\begin{itemize}
\item the vector deformation $y$ is discretized by $P^1(\mathcal{T})$ -- continuous and piecewise linear triangular elements, 
\item the scalar variable $z$ is discretized by $P_{\{0, 1 \}}^0(\mathcal{T})$ -- piecewise constant triangular elements attaining values $0$ or $1$. 
\end{itemize}
Then the time sequence of  incremental minimizations \eqref{incremental}  rewrites as:
 \begin{align}\label{discrete_problem}
  \begin{aligned}
 &\text{minimize } \\
 &\sum_{T \in \mathcal{T}}  \left( z \hat W_1(F) + (1-z) \hat W_2(F) + \beta | z - z_{k-1} | \right) |T| \\
 & \qquad \qquad \qquad +  \sum_{E \in \mathcal{E_I}} \left( |z^+ - z^-| ( \alpha_{\rm i} + \alpha_{\rm s} |\cof \mathbb{F} n | ) \right) |E|  \\
 &\text{over } y \in [P^1(\mathcal{T})]^2, z \in  P^0_{\{0,1\}} (\mathcal{T}).
   \end{aligned} 
 \end{align}
All terms  in \eqref{discrete_problem} are either constant on each triangle $T \in \mathcal{T}$ or constant on each edge $E \in \mathcal{E_I}$. Recall the deformation and surface deformation gradients $F$ and $\mathbb{F}$ are the functions of the searched deformation $y$ and $\overline \Omega \subset \mathbb{R}^2$. The dimension of the minimization problem \eqref{discrete_problem} is equal to 
$$2|\mathcal{N_F}|+|\mathcal{T}|,$$
where $|\mathcal{N_F}|$ and $|\mathcal{T}|$ denotes the number of free nodes (nodes that are not Dirichlet boundary nodes) and the number of triangles of the triangulation, respectively.
 

\begin{remark}[Simplifications] The dissipation related term simplifies as 
\begin{equation}
| z - z_{k-1} | 
= s(z_{k-1}) (z - z_{k-1}),
\end{equation} 
where
$s:\{0,1\} \rightarrow \{-1,1\}$ such that $s(0):=1$ and $s(1):=-1$. This replaces the nondifferentiable term by the linear term with the coefficient $s(z_{k-1})$ depending on the previous value $z_{k-1}$. In order to shorten the notation, we introduce abbreviations.
Here, we used an abbreviation
\begin{eqnarray}
&\hat W_{1,2}(F) := \hat W_1(F) - \hat W_2(F), \\
&\alpha: = \alpha_{\rm i} + \alpha_{\rm s} |\cof \mathbb{F} n |.
\end{eqnarray}
\end{remark}

Let us study in detail the  minimization of \eqref{discrete_problem} over $z$  for given $y$. If interfacial energy is completely neglected, $\alpha_{\rm i} = \alpha_{\rm s} = 0$, then \eqref{discrete_problem} is reduced to ($z$-independent terms are omitted): 
\begin{align}\label{discrete_problem_simplified}
  \begin{aligned}
 &\text{minimize } \\
 &\sum_{T \in \mathcal{T}}  (\hat W_{1,2}(F) + \beta \, s(z_{k-1}) ) \, z \, |T| \\
&\text{over } z \in  P^0_{\{0,1\}} (\mathcal{T}).
   \end{aligned} 
 \end{align}
This problem decouples to separate triangles  $T \in \mathcal{T}$ and has a minimizer
 \begin{align} \label{discrete_problem_z_only_alpha_0_exact_form} 
 z_{\rm{min}}=
   \begin{TC}
     \displaystyle 0 
     &\text{on triangles where $\hat W_{1,2}(F) + \beta \, s(z_{k-1})> 0$,} \\
     1 &\text{on triangles where $\hat W_{1,2}(F) + \beta \, s(z_{k-1})< 0$.   }
   \end{TC}
\end{align}
Note that if the condition $\hat W_{1,2}(F) + \beta \, s(z_{k-1}) = 0$ is satisfied, there is not a unique minimizer and both $z_{\rm{min}}=0$ and $z_{\rm{min}}=1$ are minimizers of \eqref{discrete_problem_simplified}. If the interfacial energy is taken into account, we get: 
\begin{align}\label{discrete_energy_z_only_no_multipliers}
  \begin{aligned}
 &\text{minimize } \\
 &\sum_{T \in \mathcal{T}} (\hat W_{1,2}(F) + \beta \, s(z_{k-1}) ) \, z \, |T| 
 +  \sum_{E \in \mathcal{E_I}}  
  \alpha|  z^+ - z^- | |E| \\
 &\text{over } z \in  P^0_{\{0,1\}}(\mathcal{T}).
   \end{aligned} 
 \end{align}
We introduce a new variable (Lagrange multiplier)
$$\sigma \in P^0_{\{0,1\}}(\mathcal{E_I})$$
satisfying additional constraints:
\begin{equation}\label{con1}
\begin{split}
z^+ - z^- - \sigma  &\leq 0 \\
z^- - z^+ - \sigma  &\leq 0
\end{split}
\qquad \text{on all edges }  E \in \mathcal{E_I}.
\end{equation}
Then the problem: 
 \begin{align}\label{discrete_problem_z_only}
   \begin{aligned}
& \text{minimize }  \\
&\sum_{T \in \mathcal{T}} (\hat W_{1,2}(F) + \beta \, s(z_{k-1}) ) \, z \, |T|  
+  \sum_{E \in \mathcal{E_I}}  \alpha \, \sigma |E| \\
&\text{over } z \in  P^0_{\{0,1\}}(\mathcal{T}), \sigma \in  P^0_{\{0,1\}}(\mathcal{E_I}) \\
&\text{subject to constraints } \eqref{con1}\\
     \end{aligned} 
 \end{align}
 has the minimizer part $\tilde z_{\rm{min}} \in  P^0_{\{0, 1\}}(\mathcal{T})$ equal to the minimizer of \eqref{discrete_energy_z_only_no_multipliers} and the other part $\sigma_{\rm{min}} \in P^0_{\{0,1\}}(\mathcal{E_I})$
satisfies
 \begin{align*}
\sigma_{\rm{min}} &= | z^+_{\rm{min}} -  z^-_{\rm{min}} |  & \text{on all edges }  E \in \mathcal{E_I}.
 \end{align*}
\begin{remark}
Let us consider \eqref{discrete_energy_z_only_no_multipliers} over a larger convexified set $P^0_{[0,1]}(\mathcal{E_I})$. The detailed analysis shows that the minimizer must satisfy 
$z_{\rm{min}} \in P^0_{\{0,1\}}(\mathcal{E_I})$ (none of its vector entries lies between 0 and 1). Namely,  \eqref{discrete_energy_z_only_no_multipliers} can be written as a piecewise linear  function on $[0,1]^{|\mathcal{T}|}$. Subsets of $[0,1]^{|\mathcal{T}|}$ on which the function is linear  are polygons whose extreme points coincide with extreme points of $[0,1]^{|\mathcal{T}|}$.  Therefore, it is possible to replace admissible sets 
$$P^0_{\{0,1\}}(\mathcal{T}), \quad  P^0_{\{0,1\}}(\mathcal{E_I})  $$ 
by convexified sets  
$$P^0_{[0,1]}(\mathcal{T}), \quad P^0_{[0,1]}(\mathcal{E_I}) $$ 
in practical computations. Then, the resulting linear programming problem is computationally easier to solve than the integer linear programming problem \eqref{discrete_problem_z_only}.
\end{remark}

After plugging in the implicit relations \eqref{discrete_problem_z_only} or \eqref{discrete_problem_z_only_alpha_0_exact_form} of $z_{\rm{min}}$ on $C$ (or equivalently on $F$),  it is possible to formally understand and minimize  the functional of \eqref{discrete_problem} over $y \in P^1(\mathcal{T})$ only. The functional in \eqref{discrete_problem} is convex in $C$ but known to be nonconvex in $F$  and thus nonconvex in $y$. Therefore, we obtain only a particular local minima $y$ (and the corresponding pair $(y,z)$) in our numerical computations. \\

For FEM computations, we consider a uniform triangular mesh $\mathcal{T}$ with 128 rectangles grouped in 8 horizontal layers and consisting of an underlying triangular mesh with 256 triangles. 
Deformed elements are visualized on the underlying rectangular mesh, where diagonal edges are dropped out for better readability, see Figure \ref{mesh} (right). 

The time-dependence of the stored energy is realized through the time-dependent Dirichlet boundary conditions on the deformation $y$. The boundary conditions are parametrized by a function $a: [0, T] \rightarrow \R$. We consider one period (T=16) of the triangular wave function with amplitude 1 defined as a linear interpolation function with discrete values: 
$$
\begin{tabular}{ l | c | c | c | c | c }	
  t & 0 & 4 & 8 & 12 & 16 \\
	\hline
  $a(t)$ & 0 & 1 & 0 & -1  & 0 \\
\end{tabular}
$$ 
Together, we compute 16 incremental time-step problems (N=16). Parameters of the Mooney-Rivlin material model \eqref{density_Mooney-Rivlin_C} were chosen as $\alpha=1,  \delta_1=1$ (it implies $\delta_2=4$), the structural parameter $\epsilon = 0.3$, and the dissipation coefficient $\beta=0.1$. For simplicity, we put $\alpha_{\rm{s}}=0$. In both examples, we construct the initial fields $u(0,x), z(0,x)$ at $t=0$ by solving the problem \eqref{discrete_problem} without the dissipation (the term $\beta | z - z_{k-1} |$ is omitted) and starting the minimization from the zero displacement field $u(0,x)$; the field $z(0,x)$ is prescribed as specified below. By this way, we obtain what can be called the initial relaxed microstructure.


\subsection{Example 1}
In the first (and a bit artificial) example, we study the influence of the interfacial energy term on the general response. A displacement boundary condition  is prescribed on the full domain boundary, $\Gamma_{\rm D} = \partial \Omega$, as $$ u(t,x) = 0.3 \, a(t)(x_2 - 0.5,0) \qquad \forall x=(x_1, x_2) \in \Gamma_{\rm D}.$$ (Limit positions of the domain boundary correspond to those at Figure \ref{mesh_stretching}.) We study two modeling cases with interface energy parameters: 
\begin{itemize}
\item[] $\alpha_{\rm{i}}=0$ (see Figure \ref{example1a}) 
\item[] $\alpha_{\rm{i}}=0.003$ (see Figure \ref{example1b})
\end{itemize}
At $t=0$, a (pseudo-)random distribution of variants is generated (the same for both cases), so that fraction of variants is approximately the same (less than $5\%$ difference), see Figure \ref{example1a} at $t=0$. Notice that thanks to the interfacial energy contribution, the initial relaxed microstructure differs in Figure \ref{example1b} -- more compact clusters of variants nucleate. With progressive stretching of the the domain, fraction of currently more favorable variant always increases as expected. However, the process occurs in a generally haphazard manner when interfacial term is not considered, whereas the less favorable variant disappears predominantly by shrinking of the corresponding variant-region in the other case, so that the interfacial energy contribution monotonically decreases. In both cases, hysteretic behavior is manifested, i.e. the fraction of variants is always ``delayed" with respect to the equilibrium situation with no dissipation term. This can be best observed when the domain takes the rectangular shape ($t=8,16$): the more favorable variant for preceding loading phase dominates the sample (approx. $3:1$ ratio for no interfacial energy and more than $90\%$ in the other case) with a pronounced effect of interfacial term penalizing fragmentation.

\subsection{Example 2}
The second example deals with a transition from a generic laminated structure (e.g. a martensitic laminate obtained by cooling from austenite) to a structure which is more favorable to the applied loading. 
We study one modeling case only with the interface energy parameter: 
\begin{itemize}
\item[] $\alpha_{\rm{i}}=0.001$ (see Figure \ref{example2}). 
\end{itemize}
The initial configuration consists of alternating horizontal strips of single variants, see Figure \ref{example2} at time $t=0$. The vertical movement of the right edge is prescribed, the left edge is fixed and a periodic boundary condition is applied on the top and bottom edges, so that a simple shear of an (vertically) infinite 2D SMA ribbon is mimicked. Particularly, 
$$ u(t,x) = 0.4\,a(t)(0,x_1), \qquad \forall x=(x_1, x_2) \in \Gamma_{\rm D},$$
$$ u(t,x) = u(t,x-(0,1)) \qquad \forall x=(x_1, x_2) \in {\Gamma_t},$$
where
\begin{eqnarray*}
&&\Gamma_{\rm D} = \{(x_1,x_2) | x_1 = 0\} \cup \{(x_1,x_2) | x_1 = 2\}, \\
&&\Gamma_{\rm t} = \{(x_1,x_2) | x_1 \in (0,1), x_2 = 1\}.
\end{eqnarray*}

During the first shift of the right edge upwards, compact regions of the more favorable variant ($F_1$) grow from both side edges inward, whereas the central part deforms only elastically without any change in the variant pattern, see $t=2, 4$. Let us note that this may be conceptually related to the notion of localized reorientation as suggested in \cite{Liu-2000}. When the right edge is then shifted down, first the variant boundaries do not evolve (effectively impeded by dissipation), and the structure just evolves elastically. Only when the losses from the dissipation are outweighed by the gains from further elastic relaxation and the release of the interfacial energy, one observes a growth of the currently loading-favorable variant ($F_2$) from the so-far intact central part and its spread over the whole computational domain ($t=10, 12$). During the final unloading to the initial position, the variant composition does not evolve anymore.

\begin{remark}
A Matlab code providing computational results is available for download and testing at \url{https://www.mathworks.com/matlabcentral/fileexchange/68547}.
It utilizes fast vectorized assemblies of finite elements \cite{RahmanValdman-2013} previously used in own computations with damage and plasticity models \cite{RoubicekValdman-2016, RoubicekValdman-2017} and thermo-magnetic models \cite{KruzikValdman-2018}.

\end{remark}

\def \wid {1}
\def \rate{0.68}
\def \hsp{1.0}
\def \vsp{0.5}
\begin{figure} 
\centering 
\begin{minipage}[t]{\wid\linewidth}\centering
\includegraphics[width=\rate\textwidth]{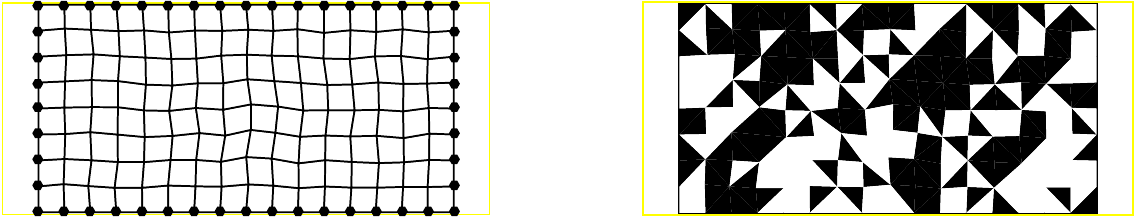}
\hspace{\hsp cm} $t=\,0$ \vspace{\vsp cm}\\
\includegraphics[width=\rate\textwidth]{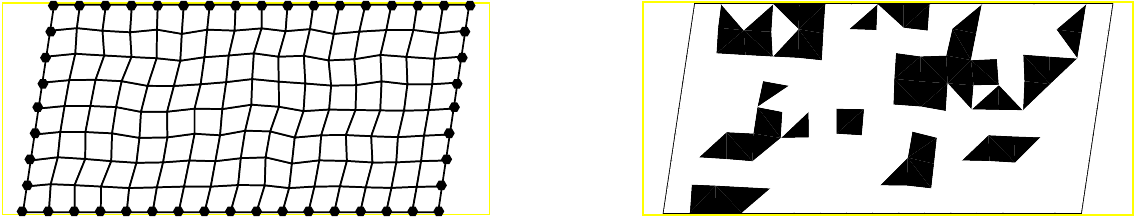}
\hspace{\hsp cm} $t=\,2$ \vspace{\vsp cm}\\
\includegraphics[width=\rate\textwidth]{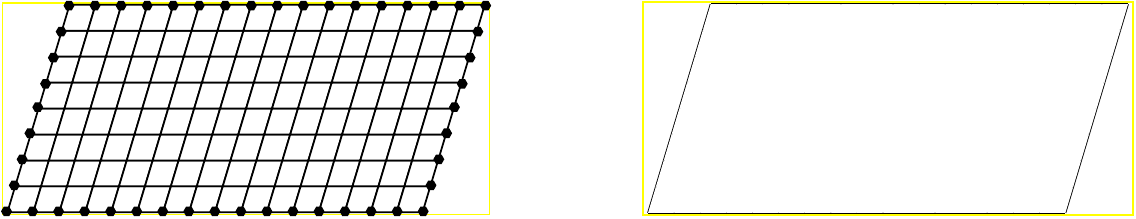}
\hspace{\hsp cm} $t=\,4$ \vspace{\vsp cm}\\
\includegraphics[width=\rate\textwidth]{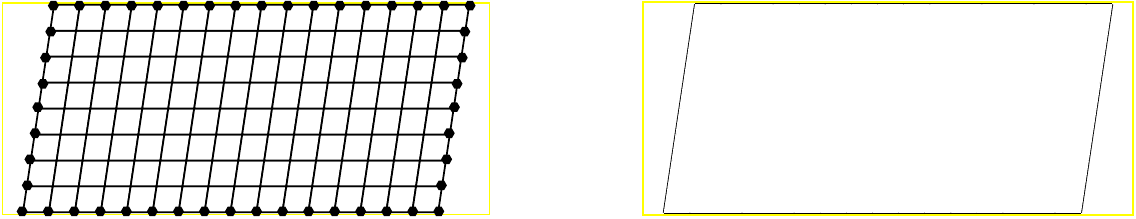}
\hspace{\hsp cm} $t=\,6$ \vspace{\vsp cm}\\
\includegraphics[width=\rate\textwidth]{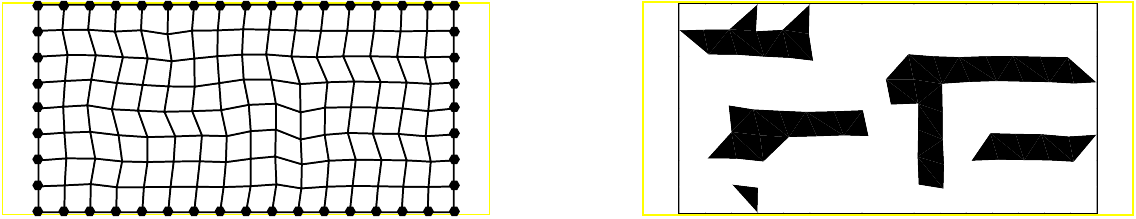}
\hspace{\hsp cm} $t=\,8$ \vspace{\vsp cm}\\
\includegraphics[width=\rate\textwidth]{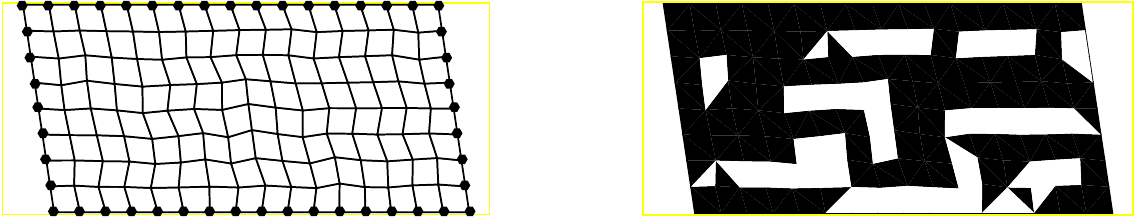}
\hspace{\hsp cm} $t=\,10$ \vspace{\vsp cm}\\
\includegraphics[width=\rate\textwidth]{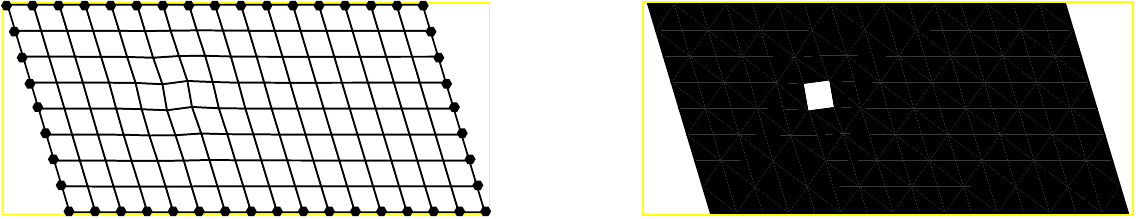}
\hspace{\hsp cm} $t=\,12$ \vspace{\vsp cm}\\
\includegraphics[width=\rate\textwidth]{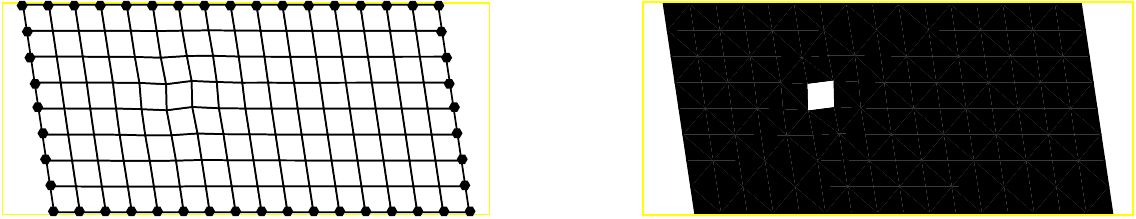}
\hspace{\hsp cm} $t=\,14$ \vspace{\vsp cm}\\
\includegraphics[width=\rate\textwidth]{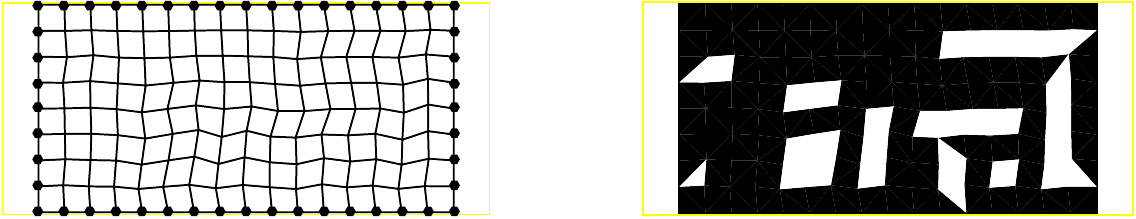}
\hspace{\hsp cm}$t=\,16$ \vspace{\vsp cm}\\
\end{minipage}
\caption{Example 1 -- $\alpha_{\rm{i}}=0$. Evolution of displacement (left) and variants distribution (right) in the domain at particular time steps.}\label{example1a}
\end{figure}

\begin{figure} 
\centering 
\begin{minipage}[t]{\wid\linewidth}\centering
\includegraphics[width=\rate\textwidth]{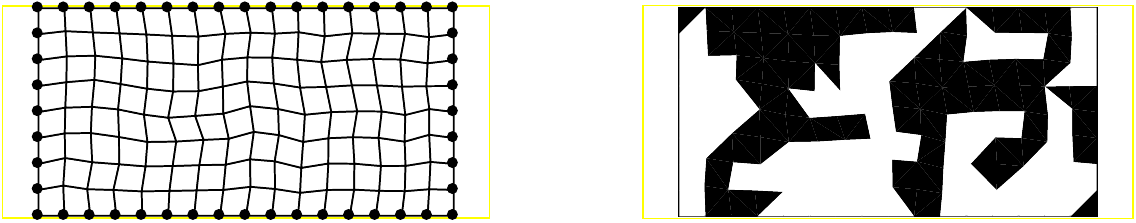}
\hspace{\hsp cm} $t=\,0$ \vspace{\vsp cm}\\
\includegraphics[width=\rate\textwidth]{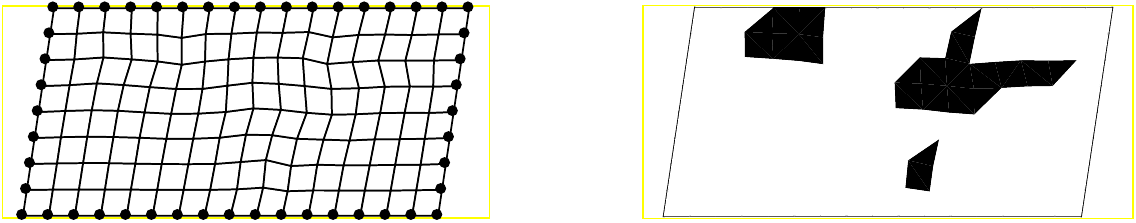}
\hspace{\hsp cm} $t=\,2$ \vspace{\vsp cm}\\
\includegraphics[width=\rate\textwidth]{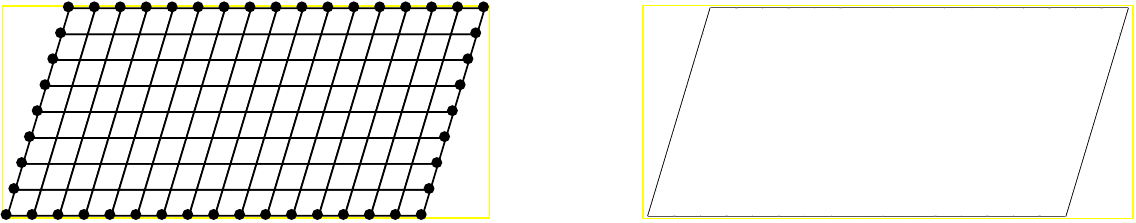}
\hspace{\hsp cm} $t=\,4$ \vspace{\vsp cm}\\
\includegraphics[width=\rate\textwidth]{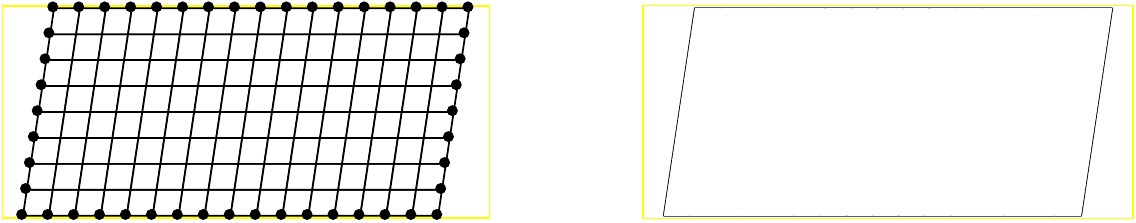}
\hspace{\hsp cm} $t=\,6$ \vspace{\vsp cm}\\
\includegraphics[width=\rate\textwidth]{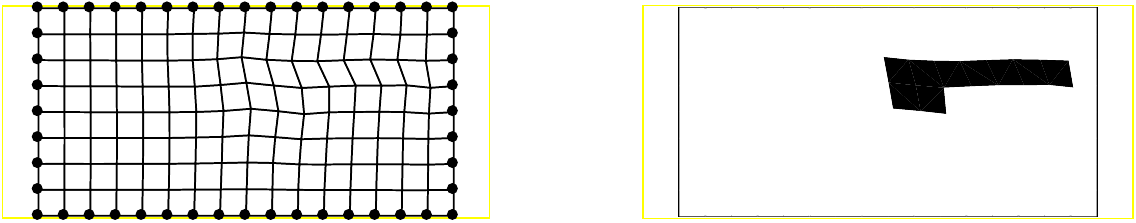}
\hspace{\hsp cm} $t=\,8$ \vspace{\vsp cm}\\
\includegraphics[width=\rate\textwidth]{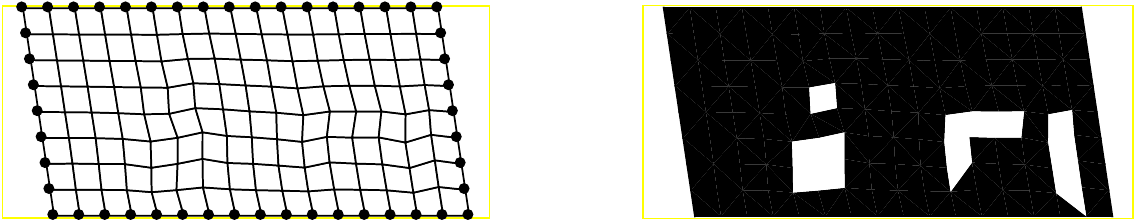}
\hspace{\hsp cm} $t=\,10$ \vspace{\vsp cm}\\
\includegraphics[width=\rate\textwidth]{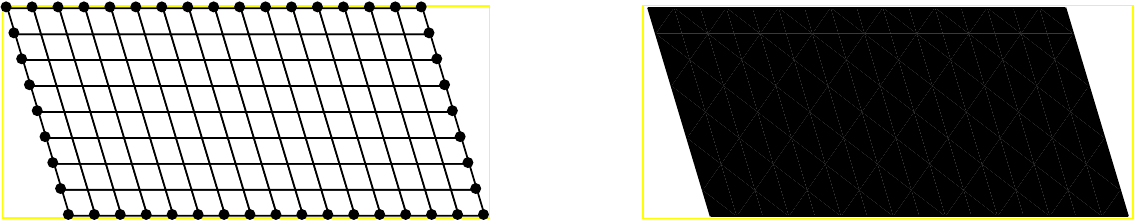}
\hspace{\hsp cm} $t=\,12$ \vspace{\vsp cm}\\
\includegraphics[width=\rate\textwidth]{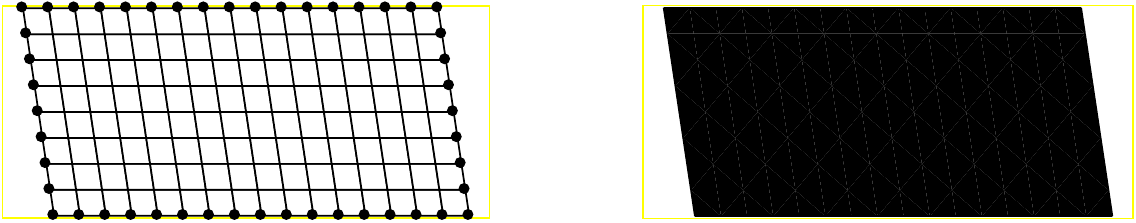}
\hspace{\hsp cm} $t=\,14$ \vspace{\vsp cm}\\
\includegraphics[width=\rate\textwidth]{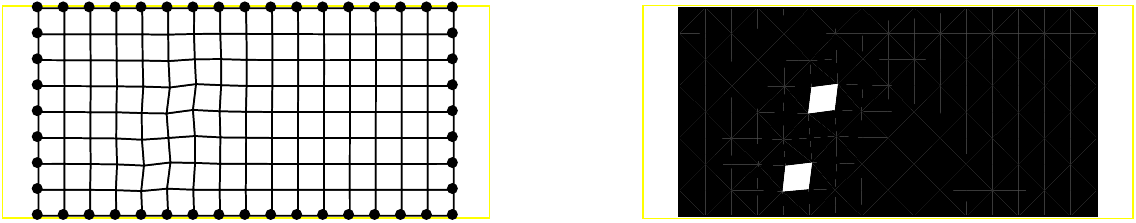}
\hspace{\hsp cm}$t=\,16$ \vspace{\vsp cm}\\
\end{minipage}
\caption{Example 1 -- $\alpha_{\rm{i}}=0.003$. Evolution of displacement (left) and variants distribution (right) in the domain at particular time steps.}\label{example1b}
\end{figure}

\def \wid {1}
\def \rate{0.5}
\def \hsp{1.0}
\def \vsp{0.35}
\begin{figure} 
\centering 
\begin{minipage}[t]{\wid\linewidth}\centering
\includegraphics[width=\rate\textwidth]{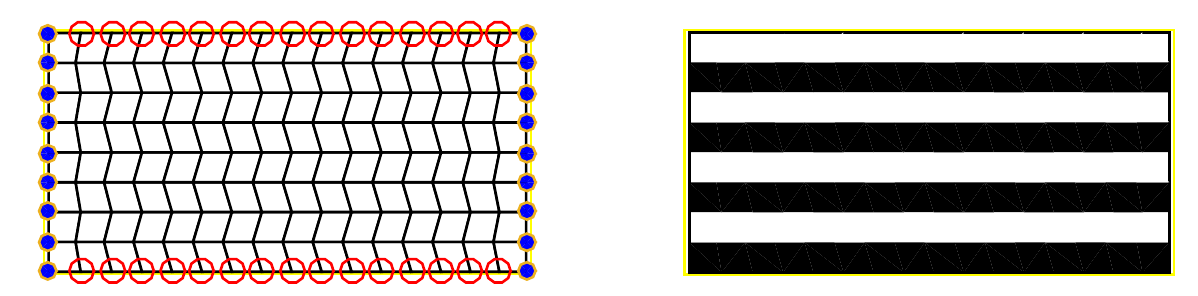}
\hspace{\hsp cm} $t=\,0$ \vspace{\vsp cm}\\
\includegraphics[width=\rate\textwidth]{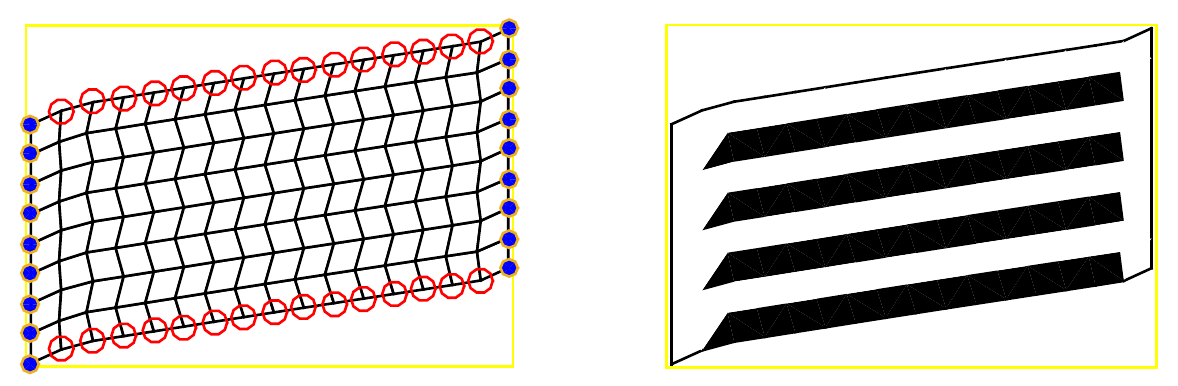}
\hspace{\hsp cm} $t=\,2$ \vspace{\vsp cm}\\
\includegraphics[width=\rate\textwidth]{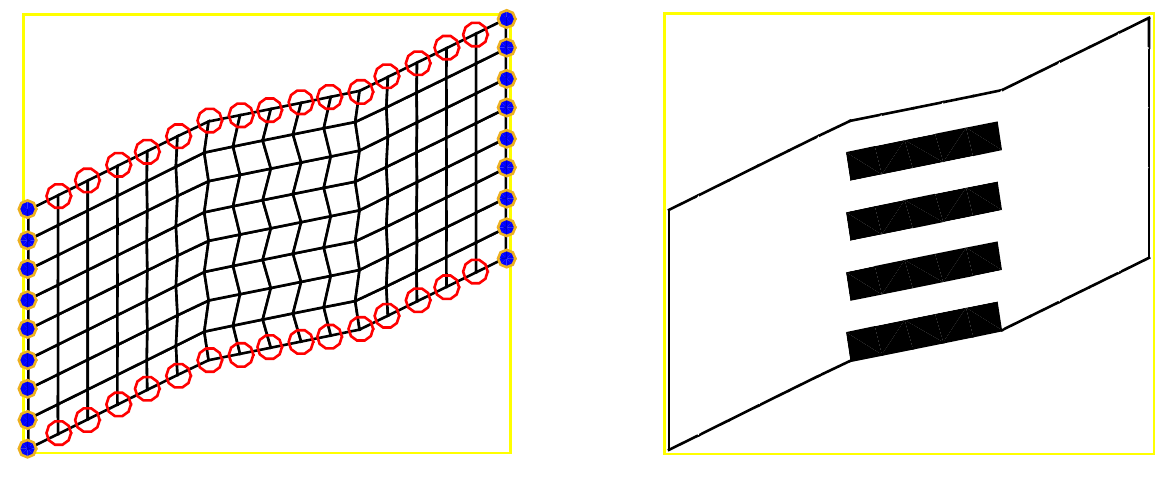}
\hspace{\hsp cm} $t=\,4$ \vspace{\vsp cm}\\
\includegraphics[width=\rate\textwidth]{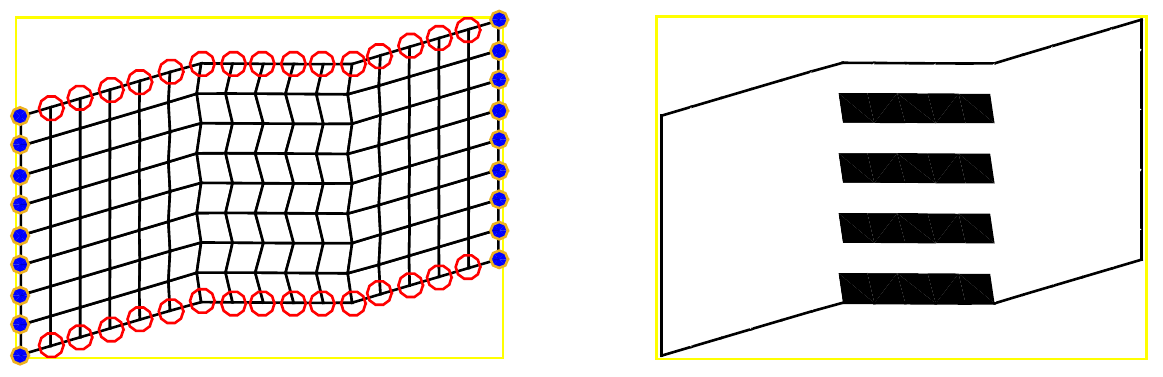}
\hspace{\hsp cm} $t=\,6$ \vspace{\vsp cm}\\
\includegraphics[width=\rate\textwidth]{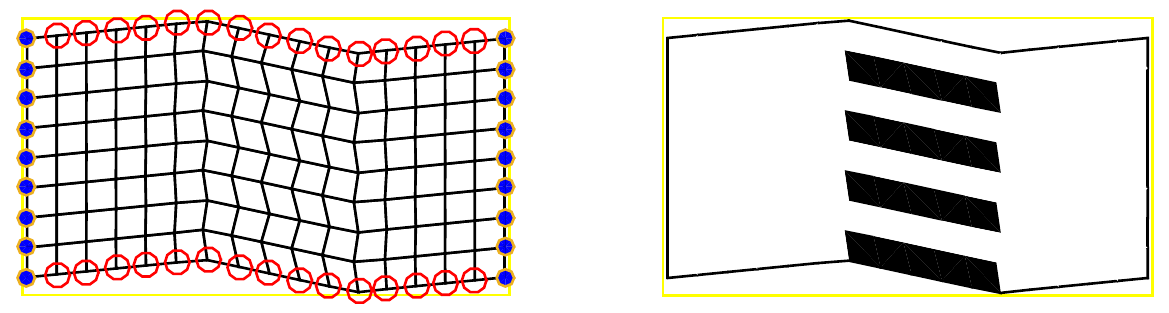}
\hspace{\hsp cm} $t=\,8$ \vspace{\vsp cm}\\
\includegraphics[width=\rate\textwidth]{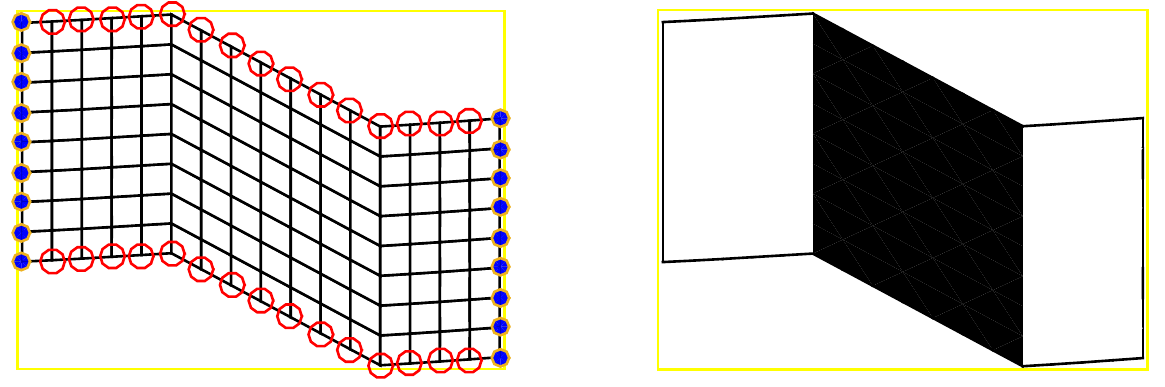}
\hspace{\hsp cm} $t=\,10$ \vspace{\vsp cm}\\
\includegraphics[width=\rate\textwidth]{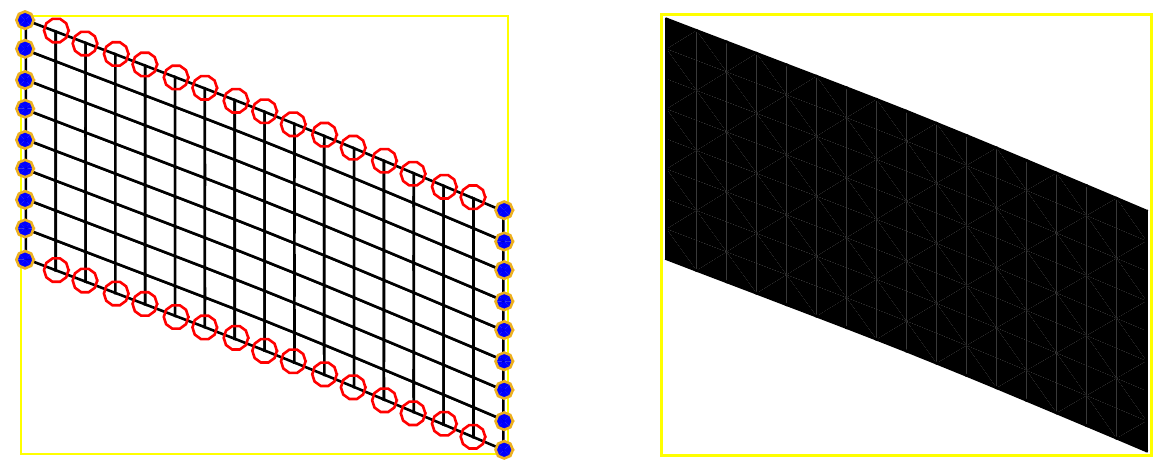}
\hspace{\hsp cm} $t=\,12$ \vspace{\vsp cm}\\
\includegraphics[width=\rate\textwidth]{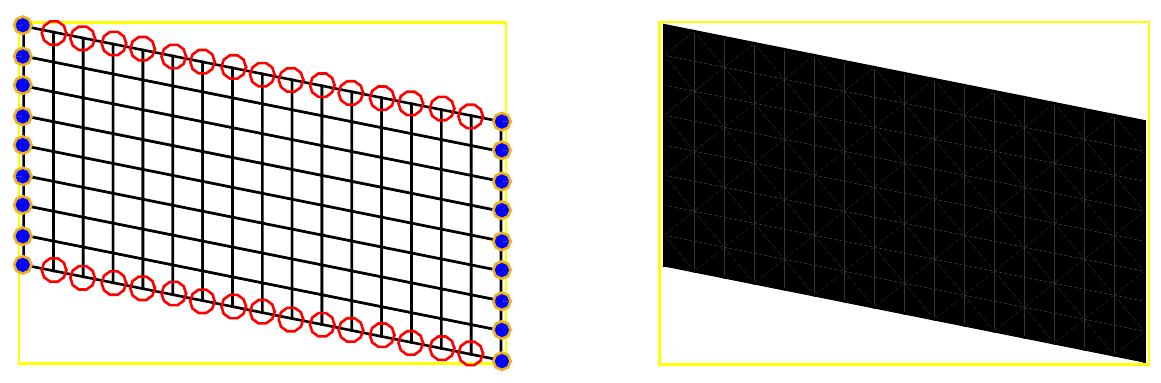}
\hspace{\hsp cm} $t=\,14$ \vspace{\vsp cm}\\
\end{minipage}
\caption{Example 2. Evolution of displacement (left) and variants distribution (right) in the domain at particular time steps. Nodes at boundaries with periodic boundary condition marked by red circles.}\label{example2}
\end{figure}


\bigskip
 {\bf Acknowledgment.}  The financial support by the Czech Science Foundation (GA\v{C}R)  through the project 18-03834S is gratefully acknowledged.

\end{document}

%% file: all-def.tex
  { \end{itemize} }
  { \end{enumerate} }

\newcommand{\qqquad}[0]{\qquad\qquad}

\newcommand{\NT}[0]{\notag}

\newcommand{\HW}[1]{} 

\newcommand{\IDEE}[1]{}
 
\allowdisplaybreaks[3]

\def\ED{\end{document}}

\newenvironment{TC} {\left \{\begin{array}{ll}} {\end{array} \right.}
\newcounter{pcounter} 
 
 \newcommand{\BS}[0]{\backslash}

\newcommand{\p}{\partial}

\newcommand{\DD}{{\mathcal D}}
\newcommand{\EE}{{\mathcal E}}

\newcommand{\HH}{{\mathcal H}}

\newcommand{\QQ}{{\mathcal Q}}

\newcommand{\VV}{{\mathcal V}}

\newcommand{\ZZ}{{\mathcal Z}}

\newcommand{\R}{\mathbb{R}}
\newcommand{\N}{\mathbb{N}}

\newcommand{\eps}{\epsilon}

\renewcommand{\phi}{\varphi}

\newcommand{\Gam}{\Gamma}

\newcommand{\Ome}{\Omega}

\DeclareMathOperator*{\cof}{cof}

\DeclareMathOperator*{\SO}{SO}

\def\XXint#1#2#3{{\setbox0=\hbox{$#1{#2#3}{\int}$}
\vcenter{\hbox{$#2#3$}}\kern-.5\wd0}}

\usepackage{color}
\definecolor{verylightblue}{rgb}{0.95, 0.95, 0.95}  
\definecolor{lightblue}{rgb}{0.7, 0.7, 1}

\definecolor{eqyellow}{rgb}{0.9375,0.8984,0.5469}
\definecolor{subeqyellow}{rgb}{1,0.9373,0.8353}

\definecolor{mygreen}{rgb}{0.3, 0.6, 0.3} 
\definecolor{verylightgreen}{rgb}{0.95, 0.95, 0.95} 
\definecolor{verydarkgreen}{rgb}{0, 0.5, 0}
\definecolor{darkgreen}{rgb}{0.85, 0.85, 0.85}  
\definecolor{mydarkgreen}{rgb}{0, 0.5, 0} 

\definecolor{mybrown}{rgb}{0.85, 0.4, 0.3}
\definecolor{verylightbrown}{rgb}{0.98, 0.72, 0.58}
\definecolor{verydarkbrown}{rgb}{0.44, 0.26, 0.26}

\definecolor{orange}{rgb}{1, 0.5, 0}

\definecolor{BurntOrange}{rgb}{1,0.356,0}

\definecolor{mydarkred}{rgb}{1,0.086,0.255}

\definecolor{RoseVYDP}{rgb}{0.84,0.086,0.255}

\definecolor{dgreen}{rgb}{0, 0.8, 0.5}     

\definecolor{CanaryBRT}{rgb}{1,0.76,0.26}

\definecolor{cyan}{rgb}{0, 1, 1}
\definecolor{verylightgray}{rgb}{0.95, 0.95, 0.95}
\definecolor{verylightgray}{rgb}{0.95, 0.95, 0.95}
\definecolor{verylightred}{rgb}{1, 0.8, 0.78}
\definecolor{verylightyellow}{rgb}{0.99, 0.98, 0.5}

\newcommand{\ignore}[1]{{}}

